\documentclass[11pt, reqno]{amsart}
  
\usepackage{lipsum}
\usepackage{a4wide}
\usepackage{amssymb} 
\usepackage{amsmath}
\usepackage{amsthm} 
\usepackage[all,cmtip]{xy}
\usepackage{graphicx}
\usepackage{tikz} 
\usepackage{amscd}
\usepackage{hyperref}
\usepackage{enumitem}
\hypersetup{colorlinks,linkcolor={blue},citecolor={blue},urlcolor={red}}
\theoremstyle{plain}
\numberwithin{equation}{section}

\newcommand{\gl}{\operatorname{GL}}

\newcommand{\tr}{\operatorname{tr}}
\newcommand{\End}{\operatorname{End}}

\newcommand{\e}{\operatorname{{\epsilon}}}
\newcommand{\tensor}{\otimes}
\newcommand{\grass}{\Lambda_\epsilon}
\newcommand{\C}{\mathbb{C}}

\newcommand{\N}{\mathbb{N}}

\newtheorem{theorem}{Theorem}[section]
 
\newtheorem{lemma}[theorem]{Lemma}
\newtheorem{remark}[theorem]{Remark}
\newtheorem{proposition}[theorem]{Proposition}

\newcommand\blfootnote[1]{%
  \begingroup
  \renewcommand\thefootnote{}\footnote{#1}%
  \addtocounter{footnote}{-1}%
  \endgroup
}
\begin{document}

\title{Mixed tensor invariants of Lie color algebra}
\author{santosha Pattanayak, Preena Samuel}

\maketitle
\begin{center}
    \small{ Department of Mathematics and Statistics, IIT Kanpur,
Kanpur, India\\
santosha@iitk.ac.in, preena@iitk.ac.in}
\end{center}

\begin{abstract}
 In this paper, we consider the mixed tensor space of a $G$-graded vector space where $G$ is a finite abelian group. We obtain a spanning set of invariants of the associated symmetric algebra under the action of a color analogue of the general linear group which we refer to as the general linear color group. As a consequence, we obtain a generating set for the polynomial invariants, under the simultaneous action of the general linear color group,  on color analogues of several copies of matrices. We show that in this special case, this is the set of trace monomials, which coincides with the set of generators given by Berele in \cite{Berele2}.
\end{abstract}

\section{Introduction}

Lie color algebras were introduced as ‘generalized Lie algebras’ in 1960 by Ree \cite{R}, being also called color Lie superalgebras (see \cite{BMPZ}). Since then, this kind of algebra has been an object of constant interest in mathematics, being also remarkable for the important role played in theoretical physics, especially in conformal field theory and supersymmetries. Lie color algebras have close relation with Lie superalgebras. Any Lie superalgebra is a Lie color algebra defined by the simplest nontrivial abelian group $\mathbb Z_2$, while any Lie color  algebra defined by a finitely generated abelian group admits a natural Lie superalgebra structure. Unlike Lie algebras and Lie superalgebras, structures and representations of Lie color algebras are far from being well understood and also there is no general classification result on simple Lie color algebras. Some recent interest related to their representation theory and related graded ring theory can be found in \cite{CSO}, \cite{Z}, \cite{SZZ}.

 In \cite{P}, Procesi studied tuples $(A_1, \cdots ,A_k)$ of endomorphisms of a finite-dimensional vector space up to simultaneous conjugation by studying the corresponding ring of invariants. He showed that for an algebraically closed field $F$ of characteristic zero the algebra of invariants $F[(A_i)_{jl}]^{GL_r(F)}$ can be generated by traces of monomials in $A_1, \cdots A_k$. 
The main tool used in the work of Procesi, in order to describe the invariants, is the Schur–Weyl duality. This tool was used in a similar way also to study more complicated algebraic structures than just a vector space equipped with endomorphisms. In \cite{dks}, Datt, Kodiyalam and Sunder applied this machinery to the study of finite-dimensional complex semisimple Hopf algebras. They were able to obtain a complete set of invariants for separating isomorphism classes of complex semi-simple Hopf algebras of a given dimension. This  they  accomplished by giving an explicit spanning set for the invariant ring of the mixed tensor space. These are called  ``picture invariants". The picture invariants obtained in the case of complex semisimple Hopf algebras were also obtained by using techniques from Geometric Invariant Theory in \cite{M}. In \cite{KR} these picture invariants were used to describe a finite collection of rational functions in the structure constants of a Lie algebra, which form a complete set of invariants for the isomorphism classes of complex semisimple Lie algebras of a given dimension. 

More recently, the invariant ring of the mixed tensor spaces was used to extend results of \cite{dks} to separate isomorphism classes of complex Hopf algebras of small dimensions \cite{PS}. The invariant ring of the mixed tensor space arises naturally in obtaining invariants for any such isomorphism classes. Here, we look at the $G$-graded analogue of this invariant ring. We define ``graded picture invariants" which $\Lambda_{\epsilon}$-linearly span this ring. 


In \cite{FM}, Fischman and Montgomery proved an analogue of the double centralizer theorem for coquasitriangular Hopf algebras and as a consequence they proved the double centralizer theorem for the general linear Lie color
algebra $\mathfrak{gl}_{\epsilon}$. In \cite{Moon}, Moon reproves the double centralizer theorem for $\mathfrak{gl}_{\epsilon}$ by relating its centralizer algebra to that of the general Lie superalgebra. For a $G$-graded vector space $V$, we set $U:=V\tensor_\C\Lambda_{\epsilon}$.
The {\it general linear color group} is defined as the group of invertible elements in the set of all degree preserving $\Lambda_{\epsilon}$-linear maps on $U$. We denote this group by $\gl_\epsilon(U)$. Using the results of \cite{FM}, Berele obtained  a graded analogue of the Schur–Weyl
duality between the general linear color group and the symmetric group in \cite{Berele2}.  This result is then used to arrive at a generating set for the $\gl_\epsilon(U)$-polynomial invariants on color analogue of multiple copies of matrices, thereby extending some of the results of \cite{P} to this setting.

In this paper, we extend Berele's results to the action of the general linear color group on the mixed tensor space $\oplus_{i=1}^sU_{b_i}^{t_i}$ where  $t_i$, $b_i$ are in $\N\cup \{0\}$ for all $i=1,\ldots s$. The case of invariants of the color analogues of $d$ copies of matrices, as described in \cite{Berele2}, may be seen as a special case of the above by taking $t_i=1=b_i$ for all $i=1\ldots,s$. We show in Theorem~\ref{t:pictures} that the ring of invariants of the $G$-graded symmetric algebra $S(\oplus_{i=1}^sU_{b_i}^{t_i})^*$, is generated by certain special invariants which are analogues of the ``picture invariants" in \cite{dks}. We then show in Theorem~\ref{t:traces} that in the special case of $t_i=1=b_i$ for all $i=1\ldots,s$, this agrees with the results of \cite{Berele2}. Viewing a  superspace as a special case of a $G$-graded space, the invariants that we obtain here coincide with those obtained in \cite[Theorem~4.1]{PS}.

  To define polynomials on a $G$-graded vector space we use the notion of $\grass$-valued polynomials over $U$, as introduced in \cite{LZ}.  This is done in \S\ref{ss:polynomials}.  This description of polynomials over $U$ is a suitable alternative to Berele's notion of a polynomial  as defined in \cite{Berele2} since we are interested in mixed tensor spaces which do not have a natural identification with matrices. In the special case of $t_i=1=b_i$ for $i=1,\ldots,s$, however, this notion of polynomials agrees with that of \cite{Berele2}. We have used this notion also in \cite{PS} to arrive at extensions of Berele's results in the supersetting. For defining this notion of a polynomial, we consider the $\grass$-module of maps to $\grass$ from the graded component of $U$ corresponding to the identity element of $G$. This module is denoted as  $\mathcal{F}(U_0,\grass)$. Then using a $G$-graded analogue of the restitution map from the $r$-fold tensor space of $U^*$ to  $\mathcal{F}(U_0, \grass)$, we call the image under this map to be the space of homogeneous polynomials of degree $r$. The $G$-graded algebra of polynomials on $U$ is then taken to be the direct sum of these spaces of homogeneous polynomials. We prove  that this algebra is isomorphic to the symmetric algebra of $U^*$, under the restitution map, analogous to {\it loc. cit.}. In this paper, we use the above notion of $\grass$-valued polynomials on $\oplus_{i=1}^sU_{b_i}^{t_i}$ and obtain a generating set for the polynomial invariants of a $G$-graded mixed tensor space.  We then show that in the special case when the mixed tensor space corresponds to several copies of the endomorphism space of $U$, the graded picture invariants are just trace monomials as given in \cite{Berele2}. This may be regarded as the $G$-graded analogue of Procesi's result in \cite{P}.

We now give an outline of the paper. In section 2 we review preliminaries of $G$-graded vector spaces and Lie color algebras and, we also recall the $G$-graded analogue of the Schur-Weyl duality. We also introduce in this section the notion of a general linear color group which is denoted as $\gl_{\epsilon}(U)$ in \cite{Berele2} and recall the Schur-Weyl duality for it. In section 3 we introduce the notion of graded picture invariants and prove that these span the space of invariants of the symmetric algebra of the dual of the mixed tensor space. Using this result we give a spanning set for the polynomial invariants of the mixed tensor space and thereby show that the trace monomials span the  polynomial invariants for the action of the general linear color group on color analogues of several copies of matrices. 

{\bf Notation:} Throughout this paper we work over the field of complex numbers $\mathbb C$. All modules and algebras are defined over $\mathbb C$ and in addition all the modules are of finite dimension. We write $\mathbb Z_2=\{\bar{0},\bar{1}\}$ and use its standard field structure. For a finite abelian group $G$, we denote the identity element of $G$ by $\circ$, while we reserve the symbols $0,1$ to denote the usual complex numbers that they represent. 
\section{Preliminaries}
\label{s:preliminaries}
\subsection{}\label{def:epsilon}{\bf Definition:} Let $G$ be a finite abelian group.  A map $\epsilon: G \times G \rightarrow \mathbb C \setminus \{0\}$ is called a skew-symmetric bicharacter on $G$ if the following identities hold, for all $f, g, h \in G$. 

(1) $\epsilon (f, g+h)=\epsilon(f,g) \epsilon(f,h)$,

(2) $\e(g+h,f)=\e(g,f)\e(h,f)$,

(3) $\e(g,h)\e(h,g)=1$.

From the definition of $\epsilon$ it follows that $\epsilon(g,\circ)=\epsilon(\circ,g)=1$ and $\epsilon(g,g)=\pm 1$, where $\circ$ denotes the identity element of $G$. 

For a bicharacter $\epsilon$ of $G$, we set $G_{\bar 0}:=\{g \in G: \epsilon(g,g)=1\}$ and $G_{\bar 1}=G \setminus G_{\bar{0}}$. Then there exists a group homomorphism $\theta:G \rightarrow \mathbb Z_2$ such that $\theta(g)=\bar 0$ for $g \in G_{\bar 0}$ and $\theta(g)=\bar 1$ for $g \in G_{\bar 1}$. Moreover, $G_{\bar 0}$ is a subgroup of $G$ of index $1$ or $2$. It easily follows that if $g \in G_{\bar 1}$ then $-g \in G_{\bar 1}$. 

{\bf Definition:} For a finite abelian group $G$, a $G$-graded vector space is a vector space $V$ together with a decomposition into a direct sum of the form $V=\oplus_{g \in G} V_g$, where each $V_{g}$ is a vector space. For a given $g \in G$ the elements of $V_{g}$ are then called homogeneous elements of degree $g$ and we write $|v|$ to denote the degree of $v$.

Any finite dimensional $G$-graded vector space 
$V=\oplus_{g \in G} V_g$ can be given a $\mathbb Z_2$-grading via $V=V_{G_{\bar 0}} \oplus V_{G_{\bar 1}}$, where $V_{G_{\bar 0}}=\oplus_{g \in G_{\bar 0}}V_g$ and $V_{G_{\bar 1}}=\oplus_{g \in G_{\bar 1}}V_g$.

{\bf Definition:} Fix a pair $(G,\epsilon)$, where $G$ is a finite abelian group and $\epsilon$ is a skew-symmetric bicharacter on $G$. A Lie color algebra $L=\oplus_{g \in G}L_g$ associated to $(G, \epsilon)$ is a $G$-graded $\mathbb C$-vector space  with a graded bilinear map $[-,-]:L \times L \rightarrow L$ satisfying 

1. $[L_g,L_h] \subset L_{g+h}$ for every $g,h \in G$,

2. $[x, y]=-\e(x, y)[y, x]$ and 

3. $\e(z, x)[x,[y, z]] + \e(x, y)[y,[z, x]] +\e(y, z)[z,[x, y]] = 0$ for all homogeneous elements $x,y,z \in L$. 

\subsection{}\label{action} Given a $G$-graded vector space $V=\oplus_{g \in G} V_g$, the $\mathbb C$-endomorphisms of $V$, $End_{\mathbb C}(V)$ is also $G$-graded, where 
$$End_{\mathbb C}(V)_g=\{f: V \rightarrow V: f(V_h) \subset V_{g+h}, \,\, for \,\, all \,\, h \in G\}$$
The general linear Lie color algebra, $\mathfrak{gl}_{\epsilon}(V)$ is defined to be $End_{\mathbb C}(V)$ with the Lie
bracket given by $[x,y]_{\epsilon}=xy-\epsilon(x,y) yx$. 

The $k$-fold tensor product $V^{\otimes k}$ of $V$ is also $G$-graded; $V^{\otimes k}=\oplus_{g \in G} (V^{\otimes k})_g$, where $(V^{\otimes k})_g=\oplus_{g=g_1+\cdots +g_k}V_{g_1} \otimes \cdots \otimes V_{g_k}$. 

We have an action of the symmetric group $S_k$ on $V^{\otimes k}$ as follows: 

The group $S_k$ is generated by the transpositions $s_1, s_2, \cdots s_{k-1}$, where $s_i=(i,i+1)$. Then 
$$s_i.(v_1 \otimes v_2 \otimes \cdots \otimes v_k)=\epsilon(g,h) (v_1 \otimes v_2 \otimes \cdots \otimes v_{i+1} \otimes v_i \otimes \cdots \otimes v_k),$$ where $v_i \in V_g$ and $v_{i+1} \in V_h$. 
More generally, for $\sigma \in S_k$ and $\underline{v}=v_1 \otimes v_2 \otimes \cdots \otimes v_k$ where each $v_i$ is a homogeneous element of $V$,
 $$ \sigma.\underline{v}=\gamma(\underline{v}, \sigma^{-1})(v_{\sigma^{-1}(1)}\otimes v_{\sigma^{-1}(2)} \otimes \cdots \otimes v_{\sigma^{-1}(k)}),$$ where $\gamma(\underline{v}, \sigma)=\prod_{(i,j) \in Inv(\sigma)}\epsilon(|v_i|,|v_j|)$, with $Inv(\sigma)=\{(i,j): i < j \,\, \text{and} \,\, \sigma(i) > \sigma(j)\}$. 
 
 We then extend the action to  $V^{\otimes k}$ by linearity. We then have the following lemma.

\begin{lemma}
    For two permutations $\sigma, \tau \in S_k$ we have $\gamma(\underline{v},\sigma\tau)=\gamma(\sigma^{-1}\underline{v},\tau)\gamma(\underline{v},\sigma)$.
\end{lemma}

\begin{proof}
Let $w_i=v_{\sigma(i)}$ for all $i$. Then $w_{\tau(i)}=v_{\sigma\tau(i)}$ for all $i$. 

The action of $\sigma$ on $\underline{v}$ is given by $\sigma.\underline{v}=\gamma(\underline{v}, \sigma^{-1})(v_{\sigma^{-1}(1)}\otimes v_{\sigma^{-1}(2)} \otimes \cdots \otimes v_{\sigma^{-1}(k)})$

We have $\tau^{-1}\sigma^{-1}.\underline{v}=\gamma(\underline{v}, \sigma\tau)(v_{\sigma\tau(1)}\otimes v_{\sigma\tau(2)} \otimes \cdots \otimes v_{\sigma\tau(k)})$.

On the other hand $\tau^{-1}\sigma^{-1}.\underline{v}=\gamma(\underline{v}, \sigma)\tau^{-1}.(v_{\sigma(1)}\otimes v_{\sigma(2)} \otimes \cdots \otimes v_{\sigma(k)})=\gamma(\underline{v}, \sigma)\tau^{-1}.\underline{w}$

$=\gamma(\underline{v}, \sigma) \gamma(\underline{w}, \tau)(w_{\tau(1)}\otimes w_{\tau(2)} \otimes \cdots \otimes w_{\tau(k)})$

$=\gamma(\underline{v}, \sigma) \gamma(\sigma^{-1}\underline{v}, \tau)(v_{\sigma\tau(1)}\otimes v_{\sigma\tau(2)} \otimes \cdots \otimes v_{\sigma\tau(k)})$

So we have the required identity. 
\end{proof}

By fixing a set of homogeneous vectors $v_1, \cdots ,v_k$ of $V$ such that $|v_i|=a_i$ for all $i$, we set $I$ to be the tuple $(a_1,\ldots, a_k)$ in $G^k$. The symmetric group $S_k$ acts on $G^k$ via  $\sigma.(a_1,\ldots, a_{k}):=(a_{\sigma^{-1}(1)},\ldots , a_{\sigma^{-1}(k)})$ for $\sigma\in S_{k}$.
  Define $\gamma(I,\sigma):=\gamma(v_{1}\tensor \cdots \tensor v_{k},\sigma)$. Then the above relation may be rephrased as, $\gamma(\sigma^{-1} I,\tau)\gamma(I,\sigma)=\gamma(I,\sigma\tau)$.

We denote by $\Phi$ the resulting homomorphism: $\Phi: \mathbb C[S_k] \rightarrow End_{\mathbb C}(V^{\otimes k})$. 

On the other hand the general linear Lie color algebra, $\mathfrak{gl}_{\epsilon}(V)$ acts on $V^{\otimes k}$ by twisted derivation: 
$$x.(v_1 \otimes v_2 \otimes \cdots \otimes v_k)=\sum_{i=1}^k (\prod_{j <i} \epsilon(\alpha, g_j)) v_1 \otimes v_2 \otimes \cdots \otimes x.v_i \otimes \cdots \otimes v_k,$$ where $x \in \mathfrak{gl}_{\epsilon}(V)_{\alpha}$ and each $v_j \in V_{g_j}$. 

We denote by $\Psi$ the resulting homomorphism: $\Psi: \mathfrak{gl}_{\epsilon}(V) \rightarrow End_{\mathbb C}(V^{\otimes k})$.

The group $G$ also act on $V^{\otimes k}$ by 
$$g(v_1 \otimes v_2 \otimes \cdots \otimes v_k)=\prod_i\epsilon(g, g_i)(v_1 \otimes v_2 \otimes \cdots \otimes v_k),$$ where $g \in G$ and each $v_i \in V_{g_i}$.

We denote by $\eta$ the resulting homomorphism: $\eta: \mathbb C[G] \rightarrow End_{\mathbb C}(V^{\otimes k})$. We then have the double centralizer theorem for the general linear Lie color algebra. 

\begin{theorem}[Fischman and Montgomery \cite{FM}] Let $A=\Phi(\mathbb C[S_k])$ and let $B$ be the subalgebra of $End_{\mathbb C}(V^{\otimes k})$ generated by $\eta(\mathbb C[G])$ and $\Psi(\mathfrak{gl}_{\epsilon}(V))$. Then $A$ and $B$ are centralizers of each other. 
\end{theorem}

\subsection{} \label{ss:generalities} Fix a pair $(G,\epsilon)$, where $G$ is a finite abelian group and $\epsilon$ is a skew-symmetric bicharacter on $G$. We define a graded algebra $\Lambda_{\epsilon}$ which generalizes the infinite Grassmann algebra in the super setting. Let $X=\cup_{g \in G}X_g$ be a $G$-graded set, where each $X_g$ is countably infinite. Let $\Lambda$ be the free $\mathbb C$-algebra generated by $X$ and we define $\Lambda_{\epsilon}:=\Lambda/I$, where $I$ is the ideal generated by the elements $xy-\epsilon(g,h)yx$, for all $g,h \in G$ and for all $x \in X_g$ and $y \in X_h$. Then $\Lambda_{\epsilon}$ is also $G$-graded. 

If $X=\{x_1, x_2, \cdots \}$, then the set $\{x_{i_1}x_{i_2} \cdots x_{i_r}: i_1 \leq i_2 \leq \cdots \leq i_r\}$ is a basis of $\Lambda_{\epsilon}$. For a fixed linear ordering of the elements of $X$, we set $\grass(N)$ to be the linear span of the basis vectors $\{x_{i_1}x_{i_2} \cdots x_{i_r}:1\leq i_1 \leq i_2 \leq \cdots \leq i_r\leq N\}$ where $N \in \mathbb Z_{\geq 0}$. We then have the following filtration: $$\Lambda_{\epsilon}(N) \subset \Lambda_{\epsilon}(N+1), \,\, \text{for all} \,\, N \geq 0.$$ 

\subsection{}\label{ss:isomorphisms}Given a $G$-graded vector space $V=\oplus_{g \in G} V_g$, we set $U:= V \otimes_{\mathbb C} \Lambda_{\epsilon}$. Then $U$ is also a $G$-graded $\Lambda_{\epsilon}$-bimodule: $U=\oplus_{g \in G}U_g$, where $U_g=\oplus_{g_1+g_2=g} V_{g_1} \otimes (\Lambda_{\epsilon})_{g_2}$; $x.u=\epsilon(g,h)ux$ where $x$ and $u$ are of degrees $g$ and $h$ respectively. 

Let $End_{\Lambda_{\epsilon}}(U):=\{T \in End_{\mathbb C}(U): T(ux)=T(u)x \,\, for \,\, u \in U, x \in \Lambda_{\epsilon} \}$. Then $End_{\Lambda_{\epsilon}}(U)$ is also $G$-graded in a natural way via $$End_{\Lambda_{\epsilon}}(U)_{g}=\{T \in End_{\Lambda_{\epsilon}}(U): T(U_h) \subseteq U_{g+h}, \,\, for \,\, all \,\, h \in G\}.$$

There is a natural embedding $V\hookrightarrow U$ given by $v\mapsto v\otimes 1$. Any $\mathbb C$-linear map $T:V \rightarrow V$ extends uniquely to an element of $End_{\Lambda_{\epsilon}}(U)$ and so we have $End_{\Lambda_{\epsilon}}(U) \cong End_{\mathbb C}(V) \otimes_{\mathbb C} \Lambda_{\epsilon}$. 

More generally,  for two $G$-graded vector spaces $V$ and $V'$, set $U:=V\tensor_{\mathbb C}\grass$ and $U':= V'\tensor_{\mathbb C} \grass$, then we denote by $Hom_{\Lambda_{\epsilon}}(U,U')$  to be the set of $\mathbb C$-linear maps from $U$ to $U'$ which commute with the right action of $\Lambda_{\epsilon}$. Then $Hom_{\mathbb C}(V, V')\tensor_{\mathbb C}\grass\cong Hom_{\grass}(U,U')$. In particular, if we denote by $U^*$ the $G$-graded space $Hom_{\Lambda_{\epsilon}}(U,\Lambda_{\epsilon})$, then $U^*\cong V^*\tensor_{\mathbb C}\grass$.

With notation as above, it may be easily seen  that $U\tensor_{\grass}U'\cong (V\tensor_{\mathbb C} V')\tensor_{\mathbb C}\grass$ via the map defined on the homogeneous elements by $v\tensor\lambda\tensor v'\tensor \lambda'\mapsto v\tensor v'\tensor \epsilon(|\lambda|,|v'|)\lambda\lambda'$. Here $|\lambda|$, $|v'|$ denote the $G$-grading of $\lambda$ and $v'$ respectively. 

There is a pairing between $U$ and $U^*$  given by $u\tensor \alpha\mapsto \epsilon(|\alpha|,|u|)\alpha(u)$ where $u\in U$ and $\alpha\in U^*$. This will be called the {\it evaluation map} and is denoted by $ev$. More generally, let $V_1,\ldots, V_k$ be $G$-graded vector spaces and let $W_i:=V_i\tensor_{\mathbb{C}}\grass$. Let $W=W_1\oplus \cdots \oplus W_k$ and $\tau$ be the permutation which takes $(1,2,\cdots, k,k+1, \cdots ,2k)$ to $(\tau(1), \tau(2), \cdots ,\tau(k),\tau(k+1), \cdots ,\tau(2k))=(1,3,5, \cdots ,2k-1, 2,4,6, \cdots ,2k)$.  Then the group $S_{2k}$ acts naturally on the $2k$-fold tensor product $(W\oplus W^*)^{\tensor 2k}$. Under this action the element $\tau\in S_{2k}$ induces an isomorphism between the subspaces $ W_1^*\tensor \cdots \tensor W_k^* \tensor W_1\tensor \cdots \tensor W_k$ and $ W_1^*\tensor W_1\tensor \cdots \tensor W_k^* \tensor W_k$ of $(W\oplus W^*)^{\tensor 2k}$. This isomorphism followed by the map $W_1^*\tensor W_1\tensor \cdots \tensor W_k^* \tensor W_k\to \grass$ given by $\alpha_1\tensor w_1\tensor\cdots\tensor \alpha_k\tensor w_k\mapsto \prod_i\alpha_i(w_i)$ is called the {\it evaluation map} and is denoted by $ev$. The non-degeneracy of the pairing $ev:W_1^*\tensor \cdots \tensor W_k^* \tensor W_1\tensor \cdots \tensor W_k\to \grass$ comes from noticing that this map is obtained by extending scalars to $\grass$ of a non-degenerate pairing over $\mathbb{C}$. The isomorphism induced by this non-degenerate pairing will be denoted by $\iota:W_1^*\tensor \cdots \tensor W_k^* \to (W_1\tensor \cdots \tensor W_k)^*$

\subsection{Symmetric algebra on a graded vector space}  \label{ss:symm} Let $G$ be a finite abelian group and $\epsilon$ is a skew-symmetric bicharacter on $G$. Let $V=\oplus_{g \in G} V_g$ be a $G$-graded vector space over $\mathbb{C}$. The tensor algebra on $V$ is $T(V)=\oplus_{k \geq 0} V^{\otimes k}$. 

The symmetric algebra on $V$ is defined to be $S(V)=T(V)/I(V)$, where $I(V)$ is the ideal of $T(V)$ generated by elements of the form $v \otimes w-\epsilon(g,h)w \otimes v$, where $v$ and $w$ are homogeneous elements of degrees $g$ and $h$ respectively. Note that the ideal $I(V)$ is both $\mathbb Z$-graded as well as $G$-graded.

We note that if we write $V=\oplus_{g \in G} V_g=V_{G_{\bar 0}} \oplus V_{G_{\bar 1}}$, where $V_{G_{\bar 0}}=\oplus_{g \in G_{\bar 0}} V_g$ and $V_{G_{\bar 1}}=\oplus_{g \in G_{\bar 1}} V_g$ then since $\epsilon(v,v)=-1$ for $v \in V_{G_{\bar 1}}$ we have $v \otimes v=0$ in $S(V)$. 

We define the $d$-th symmetric power of $V$, written $S^d(V)$ to be the image of $V^{\otimes d}$ in $S(V)$. Since $I(V)$ is both $\mathbb Z$ as well as $G$-graded, we have $S(V)=\oplus_{d \geq 0}S^d(V)$ and each $S^d(V)$ is $G$-graded. If $V_{G_{\bar 1}}=V$ then $S(V)$ is denoted as $\bigwedge(V)$ and it is called the exterior algebra of $V$.

We then have the following lemma which will be used in the proof of the main theorem. 

\begin{lemma}\label{l:symm}
    (1) Given any map $f:V \rightarrow W$, between two $G$-graded vector spaces, there is a unique map of $\mathbb C$-algebras $S(f):S(V) \rightarrow S(W)$ carrying $V$ to $W$.

    (2) For a $G$-graded vector space $V$, we have $S(V \otimes_{\mathbb C} \Lambda_{\epsilon})=S(V) \otimes_{\mathbb C} \Lambda_{\epsilon}$.

    (3) For two $G$-graded vector spaces $V$ and $W$ we have $S(V \oplus W)=S(V) \otimes S(W)$. 

    (4) If $V$ is a $G$-graded vector space and $\{v_1, v_2, \cdots ,v_n\}$ is a basis of $V$ consisting of homogeneous elements, then the set of all monomials of the form $v_1^{i_1}..v_n^{i_n}$ such that $\sum_{j=1}^ni_j=d$ and $0\leq i_j\leq 1$ for $v_{i_j}$ in $V_{G_{\bar 1}}$ form a basis of $S^d(V)$. 
\end{lemma}

\begin{proof} (1) The map $f$ induces a map from $T(V)$ to $T(W)$ carrying the ideal of relations in $T(V)$ to the ideal of relations in $T(W)$. So we have an induced map of $\mathbb C$ algebras $S(f):S(V) \rightarrow S(W)$, which is unique by the construction. 

(2) It is clear that in the tensor algebra level the assertion holds, that is, $T(V \otimes_{\mathbb C} \Lambda_{\epsilon})=T(V) \otimes_{\mathbb C} \Lambda_{\epsilon}$. The algebra $\Lambda_{\epsilon} \otimes_{\mathbb C} S(V)$ is obtained by factoring out the ideal generated by elements of the form $1 \otimes (x \otimes y -\epsilon(x,y)y \otimes x)$ from $T(V) \otimes_{\mathbb C} \Lambda_{\epsilon}$. The element $1 \otimes (x \otimes y -\epsilon(x,y)y \otimes x)$ corresponds to $(1 \otimes x) (1 \otimes y) -\epsilon(x,y)(1 \otimes y) \otimes (1 \otimes x)$ and these elements generate the ideal of relations in $T(V \otimes_{\mathbb C} \Lambda_{\epsilon})$. 

(3) The proof is straight forward. 

(4) We write $V=\oplus_{g \in G} V_g=V_{G_{\bar 0}} \oplus V_{G_{\bar 1}}$, where $V_{G_{\bar 0}}=\oplus_{g \in G_{\bar 0}} V_g$ and $V_{G_{\bar 1}}=\oplus_{g \in G_{\bar 1}} V_g$. Then $S(V)=S(V_{G_{\bar 0}}) \otimes \bigwedge(V_{G_{\bar 1}})$, where $S(V_{G_{\bar 0}})$ and $\bigwedge(V_{G_{\bar 1}})$ are the symmetric and exterior algebras on $V_{G_{\bar 0}}$ and $V_{G_{\bar 1}}$ respectively. We then have the required result. 
\end{proof}
\begin{remark}
   \textup{In the above lemma, for $U=V \otimes_{\mathbb C} \Lambda_{\epsilon}$, $S(U)$ is defined as the quotient of $T(U)$ by the ideal generated by the elements of the form $v\tensor w-\varepsilon(g,h)w\tensor v$ regarded as elements of $U\tensor_{\grass} U$. Further, in view of (2) above, all the other statements of the lemma also hold after extending scalars to $\grass$.}
\end{remark}
\subsection{} For $V=\oplus_{g \in G} V_g$ and $U=V \otimes_{\mathbb C} \Lambda_{\epsilon}$, we set $$M_{\epsilon}(U):=\{T: U \rightarrow U: T(U_g) \subseteq U_{g}, \,\, for \,\, all \,\, g \in G\}.$$ With the identification,  $End_{\Lambda_{\epsilon}}(U) \cong U \otimes_{\Lambda_{\epsilon}} U^*$, $M_{\epsilon}(U)$ is spanned by all $u \otimes f$, where $u$ and $f$ are homogeneous of opposite degree, i.e., $|u|=-|f|$. We define the trace function from  $M_{\epsilon}(U)$ to $(\Lambda_{\epsilon})_{G_0}$ by $tr(u \otimes f)=\epsilon(|u|, |f|)f(u)$ with the above identification. 

By choosing a basis $\{v_1, \cdots ,v_r, v_{r+1} \cdots ,v_n\}$ of $V$, where the first $r$ of them are in $V_{G_{\bar 0}}$ and the last $n-r$ are in $V_{G_{\bar 1}}$ and identifying $M_{\epsilon}(U)$ with the space of matrices we get that $$tr(A)=\sum_{i=1}^ra_{ii}-\sum_{i=r+1}^na_{ii}.$$ It satisfies $tr(AB)=tr(BA)$ for $A,B \in M_{\epsilon}(U)$ (see Lemma 3.7 of \cite{Berele2}).

{\bf Definition:} The general linear color group $GL_{\epsilon}(U)$ is defined to be the group of invertible elements in $M_{\epsilon}(U)$.

The group $GL_{\epsilon}(U)$ acts on the $k$-fold tensor product $U^{\otimes k}$ diagonally: 

$$T.(u_1 \otimes \cdots \otimes u_k)=(Tu_1 \otimes \cdots \otimes Tu_k).$$

We denote by $\Psi$ the resulting homomorphism: $\Psi:  GL_{\epsilon}(U) \rightarrow End_{\Lambda_{\epsilon}}(U^{\otimes k})$.

The action of the symmetric group on $V^{\otimes k}$ defined in \ref{ss:generalities} can be extended to an $\Lambda_{\epsilon}$-linear action on $U^{\otimes k}$ as follows: 

For $\sigma \in S_k$ and $\underline{u}=u_1 \otimes u_2 \otimes \cdots \otimes u_k \in U^{\otimes k}$ where each $u_i$ is a homogeneous element of $U$,
 $$ \sigma.\underline{u}=\gamma(\underline{u}, \sigma^{-1})(u_{\sigma^{-1}(1)}\otimes u_{\sigma^{-1}(2)} \otimes \cdots \otimes u_{\sigma^{-1}(k)}),$$ where $\gamma(\underline{u}, \sigma)=\prod_{(i,j) \in Inv(\sigma)}\epsilon(|u_i|,|u_j|)$, with $Inv(\sigma)=\{(i,j): i < j \,\, \text{and} \,\, \sigma(i) > \sigma(j)\}$.

As $U=V \otimes \Lambda_{\epsilon}$ is a $\Lambda_{\epsilon}$-bimodule, $U^{\otimes k}$ is also a $\Lambda_{\epsilon}$-bimodule. The action of $S_k$ and $\Lambda_{\epsilon}$ on $U^{\otimes k}$ commute with each other. So the $S_k$ action on $U^{\otimes k}$ extends to an action of the group algebra $\Lambda_{\epsilon}[S_k]$. 

We denote by $\Phi$ the resulting homomorphism: $\Phi:  \Lambda_{\epsilon}[S_k] \rightarrow End_{\Lambda_{\epsilon}}(U^{\otimes k})$.

Then in \cite{Berele2}, Berele proved a version of Schur-Weyl duality for the general linear color group.

\begin{theorem}[Berele]\label{t:berele}
    Let $A$ be the subalgebra of $End_{\Lambda_{\epsilon}}(U^{\otimes k})$ generated by $\Psi(GL_{\epsilon}(U))$. Then the centralizer of $A$ in $End_{\Lambda_{\epsilon}}(U^{\otimes k})$ is $\Phi(\Lambda_{\epsilon}[S_k])$.
\end{theorem}

\section{Mixed tensor spaces}\label{mixed} For a $G$-graded vector space $V$ over $\mathbb C$, one can define the mixed tensor space as the direct sum of the form $\oplus_{i=1}^s(V^{\tensor b_i}\tensor {V^*}^{\otimes t_i})$ for an $s\in \mathbb N$ and $t_i,b_i\in \mathbb N\cup \{0\}$. For simplicity of notation we denote $V^{\tensor b_i}\tensor {V^*}^{\otimes t_i}$ by $V_{b_i}^{t_i}$. Since each summand in the mixed tensor space has a $G$-grading, there is a natural $G$-grading inherited by their direct sum $\oplus_{i=1}^s(V_{b_i}^{t_i})$. Taking $U=V\tensor_\mathbb C\grass$ and $U^*=Hom_{\grass}(U,\grass)\cong V^*\tensor_{\mathbb C}\grass$, the mixed tensor space over $\grass$ is the $G$-graded space $\oplus_{i=1}^s(U^{\tensor b_i}\tensor _{\grass}{U^*}^{\otimes t_i})$ for an $s\in \mathbb N$ and $t_i,b_i\in \mathbb N\cup \{0\}$. It may be seen that $\left(\oplus_{i=1}^s(V_{b_i}^{t_i})\right)\tensor_{\mathbb C} \grass\cong \oplus_{i=1}^s(U^{\tensor b_i} \tensor_{\grass}  {U^*}^{\otimes t_i})$. We shall denote this mixed tensor space over $\grass$ by $W$. 

Let $m=\dim V_{G_{\bar 0}}$ and $n=\dim V_{G_{\bar 1}}$. Fixing an ordering for the elements of $G_{\bar 0}$ and $G_{\bar 1}$, we list the elements of $G$ as $\{\circ,g_1,g_2,\ldots\}$ with the elements in $G_{\bar 0}$ appearing first. We then arrange the basis vectors of $V$, $\{e_i\}_{i=1}^{m+n}$, such that the first $m$ vectors are from $V_{G_{\bar 0}}$ and the rest are from $V_{G_{\bar 1}}$; further, ordered linearly so that $i<j$ implies $|e_i|$ appears before $|e_j|$ under the fixed ordering of elements of $G$. Here we use the notation $|v|=h$ for $v\in V_h$. Let $\{e_i^*\}_{i=1}^{m+n}$ be the dual basis corresponding to $\{e_i\}_{i=1}^{m+n}$. We denote the image in $U$ and $U^*$ of the above basis elements under  the embedding $V\hookrightarrow U$ (and $V^*\hookrightarrow U^*$ respectively) also by the same notation. The mixed tensor space $W$ then has a basis given in terms of the above bases of $U $ and $U^*$.
Denote the element dual to the basis vector $e_{l_1}\tensor \cdots \tensor e_{l_{b_i}}\tensor e^*_{u_1}\tensor \cdots \tensor e^*_{u_{t_i}}\in U_{b_i}^{t_i}$ by $T(i)_{l_1\ldots l_{b_i}}^{u_1\ldots u_{t_i}}$. We denote the corresponding element in $W^*$ also by the same notation. Notice that $T(i)_{l_1\ldots l_{b_i}}^{u_1\ldots u_{t_i}}$ defines a linear map on $W$ via the projection $p_i:W\to U_{b_i}^{t_i}$, {\it i.e.,} $T(i)_{l_1\ldots l_{b_i}}^{u_1\ldots u_{t_i}}(v_1,\ldots,v_s)=T(i)_{l_1\ldots l_{b_i}}^{u_1\ldots u_{t_i}}(v_i)$ for $(v_1,\ldots,v_s)\in W$.  The $G$-grading of the element $T(i)_{l_1\ldots l_{b_i}}^{u_1\ldots u_{t_i}}\in W^*$ is given by $\sum_{i=1}^{b_i}h_{l_i}-\sum_{i=1}^{t_i}h_{u_i}$ where $e_{l_j}\in V_{h_{l_j}}$ and $e^*_{u_j}\in V_{h_{u_j}}^* $. The set of all $T(i)_{l_1\ldots l_{b_i}}^{u_1\ldots u_{t_i}}$, where $i=1,\ldots , s$ and $u_j$, $l_j$ are from $\{1,\cdots, m+n\}$, forms a $\grass$-basis of $W^*$. 

\subsection{Symmetric algebra of the mixed tensor space} \label{ss:mixedsymm}  Let  $S(W^*)$ be the symmetric algebra of $W^*$. We denote by $\varpi:T(W^*)\to S(W^*)$  the natural map from the tensor algebra of $W^*$ to  $S(W^*)$. We know that $S(W^*)$ has a $\mathbb{Z}$-grading given by  $\oplus_{r\in\N\cup \{0\}}S^r(W^*)$ where  $S^r(W^*)$ is the image  under $\varpi$ of $T^r(W^*)$. The restriction of $\varpi$ to $T^r(W^*)$ is denoted as $\varpi_r$. By Lemma~\ref{l:symm}(2), $S(W^*)$ can be identified with $S(\oplus_{i=1}^s(V_{b_i}^{t_i})^*)\tensor_\C \grass$; further by Lemma~\ref{l:symm}(4)  and the remarks in \S\ref{ss:generalities}, this  in turn is identified with $[S((\oplus_{i=1}^s(V_{b_i}^{t_i})^*)_{G_{\bar 0}})\tensor \bigwedge((\oplus_{i=1}^s(V_{b_i}^{t_i})^*)_{G_{\bar 1}})] \tensor_\C \grass$. Using the relations among the $T(i)_{l_1\ldots l_{b_i}}^{u_1\ldots u_{t_i}}$, $i=1,\ldots , s$ and $u_j$, $l_j\in\{1,\cdots, m+n\}$ which are given by the ideal $I(W^*)$ of $T(W^*)$, the latter identification yields a $\grass$-basis for $S(W^*)$ consisting of monomials in $T(i)_{l_1\ldots l_{b_i}}^{u_1\ldots u_{t_i}}$, $i=1,\ldots , s$ and $u_j$, $l_j$ are from $\{1,\cdots, m+n\}$ where the variables $T(i)_{l_1\ldots l_{b_i}}^{u_1\ldots u_{t_i}}$  are  arranged in order such that basis vectors of $(\oplus_{i=1}^s{V_{b_i}^{t_i}}^*)_{g_l}$ come before those of $(\oplus_{i=1}^s{V_{b_i}^{t_i}}^*)_{g_k}$ whenever $l<k$ and among them arranged from $i=1\ldots ,s $;  the degree of each variable  $T(i)_{l_1\ldots l_{b_i}}^{u_1\ldots u_{t_i}}\in (\oplus_{i=1}^s{V_{b_i}^{t_i}}^*)_{G_{\bar 1}}$ in such a monomial being either $0$ or $1$. The monomials among the above whose total degree is $r$, form a basis  for $S^r(W^*)$ over $\grass$.

For  each $r\in\mathbb N$ and a $s$-tuple $(m_1,\ldots, m_s)$ such that 
$m_1+\cdots+ m_s=r$, the tensor space $T^{m_1}({U_{b_1}^{t_1}}^*)\tensor\cdots\tensor T^{m_s}({U_{b_s}^{t_s}}^*) $ is realised as a subspace of $T^r(W^*)$. The image of the restriction of $\varpi_r$ to this subspace is denoted as $S^{m_1,\ldots,m_s}(W^*)$. If the multidegree of the element $T(i)_{l_1\ldots l_{b_i}}^{u_1\ldots u_{t_i}}\in S(W^*)$ is set to be the $s$-tuple $(0,\ldots,0,1,0\ldots,0)$ with $1$ in exactly the $i$-th position then the monomials in the above listed basis of 
$S^r(W^*)$ with multidegree $(m_1,\ldots,m_s)$ forms a basis of $S^{m_1,\ldots, m_s}(W^*)$. The space $\tensor_{i=1}^sS^{m_i}((U_{b_i}^{t_i})^*)$ can be  identified with $S^{m_1,\ldots, m_s}(W^*)$ under the map $\phi_1\tensor \cdots\tensor \phi_s\mapsto \phi_1\cdots \phi_s.$  It is easy to see that the following diagram commutes, by checking it does on the basis elements of $T^{m_1}({U_{b_1}^{t_1}}^*)\tensor\cdots\tensor T^{m_s}({U_{b_s}^{t_s}}^*)$:
\[\xymatrix{\tensor_{i=1}^s((U_{b_i}^{t_i})^*)^{\tensor m_i}\ar@{->>}[r]^{\varpi_{m_1}\tensor\ldots\tensor \varpi_{m_s}}\ar@{^{(}->}[dd]&\tensor_{i=1}^sS^{m_i}((U_{b_i}^{t_i})^*)\ar[d]^{\cong}\\&S^{(m_1,\ldots, m_s)}(\oplus_{i=1}^s(U_{b_i}^{t_i})^*)\ar@{^{(}->}[d]\\T^r(W^*)\ar@{->>}[r]^{\varpi_m}&S^r(W^*)
}\]

Note that in the above diagram, all the maps are $\gl_{\epsilon}(U)$-equivariant; indeed, since the $\gl_{\epsilon}(U)$-action on ${(U_{b_i}^{t_i})^*}^{\tensor m_i}$ (resp., ${W^*}^{\tensor r}$) commutes with the $S_{m_i}$-action (resp., $S_r$-action) and when $M=W$ or $M=U^{t_i}_{b_i}$, we note that $S^m({M^*})$ can be identified under $\varpi_{r}$ with the $\gl_{\epsilon}(U)$-invariant subspace $e(r)\left({M^*}^{\tensor r}\right)$ where $e(r):=\frac{1}{r!}\sum_{\sigma\in S_{r}}\sigma$. Thus the $\gl_{\epsilon}(U)$-action descends to $S^{r}(M^*)$ thereby making the  map $\varpi_{r}$ a $\gl_{\epsilon}(U)$-equivariant map with respect to this action.
\subsection{Main Result}\label{ss:pic}
Keeping notations as above, let $\sum_i m_ib_i=N$ and $\sum_i m_it_i=N'$. We have the following sequence of $GL_\epsilon(U)$-invariant isomorphisms:
\begin{equation}\label{eq:iso}
\left(U^{\tensor (\sum_i m_ib_i)}\tensor (U^*)^{\tensor(\sum_i m_it_i)}\right)^*\rightarrow(\tensor_{i=1}^s(U_{b_i}^{t_i})^{\tensor m_i})^*\rightarrow\tensor_{i=1}^s((U_{b_i}^{t_i})^*)^{\tensor m_i}
\end{equation}
The first isomorphism is induced on the duals by the permutation action of $\mu\in S_{N+N'}$ on the subspace$(U\oplus U^*)^{\tensor (N+N')}$ of $(U\oplus U^*)^{\tensor (N+N')}$ where $\mu$ takes $(1,\ldots, N+N')$ to $(1,\ldots b_1, N+1, \ldots , N+t_1, b_1+1, \ldots 2b_1, \cdots  )$. The isomorphism above is $\gl_{\epsilon}(U)$-equivariant since the symmetric group action on $(U\oplus U^*)^{\tensor  (N+N')}$ commutes with the $\gl_{\epsilon}(U)$-action induced on it. 
Let $k=\sum_i m_i$ and $W=W_1\oplus \cdots \oplus W_k$, where $W_i=U_{b_i}^{t_i}$. Then the group $S_{2k}$ acts naturally on the $2k$-fold tensor product $(W\oplus W^*)^{\tensor 2k}$.  The second isomorphism in the above sequence is the inverse of $\iota$ as described in \S\ref{ss:isomorphisms} between $W_1^*\tensor \cdots \tensor W_k^* $ and $(W_1\tensor \cdots \tensor W_k)^*$. More explicitly, the isomorphism is given on the dual basis by $$T_{w_1\tensor \cdots \tensor w_{\sum_i m_i}}\mapsto {p_\epsilon(\underline w,\underline w)}T(1)_{w_1}\tensor\cdots \tensor T(s)_{ w_{\sum_i m_i}}$$ where $w_i\in W_i$ is a basis vector and $T(i)_{w_i}$ is the dual vector in $W_i^*$; here for $w_i=e_{l_1}\tensor \cdots \tensor e_{l_{b_i}}\tensor e^*_{u_1}\tensor \cdots \tensor e^*_{u_{t_i}}\in U_{b_i}^{t_i}$ the notation $T(i)_{w_i}$ represents the linear map  $T(i)_{l_1\ldots l_{b_i}}^{u_1\ldots u_{t_i}}$ by our earlier notation. The value of $p_\epsilon(\underline w,\underline w)$ for $\underline w=(w_1,\ldots, w_k)$ where $w_i\in V_{|w_i|}$ is given by $\prod_{i=1}^k(\prod_{i\leq j\leq k}\epsilon(|w_i|,|w_j|))$. This isomorphism is $\gl_{\epsilon}(U)$-equivariant since the symmetric group action on $(W\oplus W^*)^{\tensor  2k}$  commutes with the $\gl_{\epsilon}(U)$-action induced on it.

 By a standard argument, the space $U^{\tensor N}\tensor {U^*}^{\tensor N'}$ has $\gl_\epsilon(U)$-invariants if and only if $N=N'$, and so does its dual, $(U^{\tensor N}\tensor {U^*}^{\tensor N'})^*$. When $N=N'$, the latter space can be identified with $\End_{\grass
}(U^{\tensor N})$ via the non-degenerate pairing $(\End_{\grass}(U^{\tensor N}))\tensor (U^{\tensor N}\tensor {U^*}^{\tensor N})\to \Lambda_{\epsilon}$ given by $\langle A,\underline{v}\tensor \underline{f}^*\rangle:= ev(A(\underline{v})\tensor \underline{f}^*)$. The latter is $ev(\nu.\tau. (A(\underline{v})\tensor \underline{f}^*))$ where $\tau\in S_{2N}$ is as in \S\ref{ss:isomorphisms} (with $k$ replaced by $N$) and $\nu\in S_{2N}$ is the permutation $(1\, 2)(3\, 4)\cdots (2N-1\ \ 2N)$. This gives a $\gl_{\epsilon}(U)$-invariant isomorphism $$\Theta:\End_{\grass}(U^{\tensor N})\to (U^{\tensor N}\tensor {U^*}^{\tensor N})^*.$$

\subsubsection{Graded picture invariants for colored spaces}
Given an $s$-tuple $(m_1,\ldots,m_s)$ of non-negative integers such that $\sum_{k=1}^s m_kt_k=\sum_{k=1}^s m_k b_k=N$ and a $\sigma\in S_N$ we associate the polynomial $\varphi_{\sigma}\in S(W^*)$ given by \begin{equation}\label{eq:picture}
   \sum_{r_1,\ldots, r_N\in \{1,\ldots, n\}}{p_\sigma(r_1,\ldots, r_N, m_1,\ldots, m_s)}\prod_{i=1}^{s}\prod_{j=1}^{m_i}T(i)_{r'_{(\sum_{p<i}m_pb_p+(j-1)b_i+1)}\cdots r'_{(\sum_{p<i}m_pb_p+jb_i)}}^{{r_{(\sum_{p<i}m_pt_p+(j-1)t_i+1})}\cdots r_{(\sum_{p<i}m_pt_p+jt_i)}} 
\end{equation}where $r'_j:=r_{\sigma j}$ and $p_\sigma(I,M)$ for an $N$-tuple $I=(r_1,\ldots,r_N)$ of $N$ elements from $\{1,\ldots  m+n\}$ and $M=(m_1,\ldots, m_s)\in (\N\cup \{0\})^s$ such that $\sum_k m_{k}t_k=\sum_k m_k b_k=N$ takes the value in $G$ given by $$\gamma(\mu^{-1}.(I,\sigma I^*),(\nu\tau\hat{\sigma}\mu)^{-1})p(\underline{w},\underline{w})$$where $(r_1,\ldots, r_N)^{*}:=(r_1^*,\ldots, r_N^*)$, indicating that these are the indexes corresponding to the basis vectors $e_{r_1}^*,\ldots,e_{r_N}^*$  in the dual space $U^*$, $\mu, \tau,\nu$ are as defined above and for a $\sigma\in S_N$ we define $\hat{\sigma}\in S_{2N}$ as 
\begin{align*}
        \hat{\sigma}(i)&=\sigma(i) \quad\textup{ for }i\leq N \\ \hat{\sigma}(i)&=i\quad\textup{ for }i>N
\end{align*} The vector $\underline{w}:=w_1 \otimes w_2 \otimes \cdots \otimes w_{\sum_{i=1}^sm_i}$ where each $w_{\sum_{p<i}m_p+j}:=e_{r_{\sum_{p<i}m_pb_p+(j-1)b_i+1}}\tensor\cdots\tensor e_{r_{\sum_{p<i}m_pb_p+jb_i}}\tensor e_{r_{\sigma^{-1}(\sum_{p<i}m_pt_p+(j-1)t_i+1)}}^*\tensor\cdots\tensor e_{r_{\sigma^{-1}(\sum_{p<i}m_pt_p+jt_i)}}^*$ for $i=1,\ldots,s; ~j=1,\ldots ,m_i$.

\noindent
The polynomials $\varphi_\sigma$ as defined above\footnote{In (\ref{eq:picture}), the monomials should be considered modulo the anti-commutativity relations in $S(W^*)$. See Remark~\ref{r:monomialterms} above.} are called the {\it graded picture invariants}.

\begin{remark}\label{r:monomialterms} \textup{We say that a variable $$T(i)_{r_{\sigma(\sum_{p<i}m_pb_p+(j-1)b_i+1)}\cdots r_{\sigma(\sum_{p<i}m_pb_p+jb_i)}}^{r_{(\sum_{p<i}m_pt_p+(j-1)t_i+1)}\cdots r_{(\sum_{p<i}m_pt_p+jt_i)}}\in S(W^*)$$ is {\it even} or {\it odd} depending on whether the degree of the variable is in $G_{\bar{0}}$ or $G_{\bar{1}}$.
    With this terminology, we note that in the above formula  the monomials with repeated odd degree variables will be identically $0$ in $S(W^*)$. 
    In particular, those indices ${r_1,\ldots, r_N}$ in the sum in (\ref{eq:picture}) are to be dropped for which there is a $1\leq i\leq s$ and $1\leq j<j'\leq m_i$ such that the above variable corresponding to these values of $i,j,j'$ turns out to be of odd degree and, $$(r_{\sigma(\sum_{p<i}m_pb_p+(j-1)b_i+1)},\ldots, r_{\sigma(\sum_{p<i}m_pb_p+jb_i)})=(r_{\sigma(\sum_{p<i}m_pb_p+(j'-1)b_i+1)},\ldots, r_{\sigma(\sum_{p<i}m_pb_p+j'b_i)}),$$ $$(r_{\sum_{p<i}m_pt_p+(j-1)t_i+1},\ldots, r_{\sum_{p<i}m_pt_p+jt_i})=(r_{\sum_{p<i}m_pt_p+(j'-1)t_i+1},\ldots, r_{\sum_{p<i}m_pt_p+j't_i}).$$}
\end{remark}
\begin{remark}
    \textup{We retain the terminology `graded picture invariants' as used in \cite{PS} since they arise from certain combinatorial diagrams, called `pictures'. This is illustrated in \cite{PS}.}
\end{remark}
\begin{theorem}\label{t:pictures}
    With the above notation, the elements $\varphi_\sigma $ linearly span $S(W^*)^{\gl_{\epsilon}(U)}$.
\end{theorem}
\begin{proof}
     As was seen earlier in this section, the space $(U^{\tensor N}\tensor {U^*}^{\tensor N'})^*$ has $\gl_{\epsilon}(U)$-invariants if and only if $N=N'$. By Theorem~\ref{t:berele}, we know that the $\gl_{\epsilon}(U)$-invariants of $\End_{\grass}(U^{\tensor N})$ are spanned over $\grass$ by $S_N$. So, via the isomorphism $\Theta$ we get invariants on $(U^{\tensor N}\tensor {U^*}^{\tensor N})^*$. Let $(m_1,\ldots,m_s)$ be an $s$-tuple of non-negative integers such that $\sum_im_it_i=\sum m_ib_i=N$, and $\sigma\in S_N$. Going by the isomorphisms in (\ref{eq:iso}) and projecting $\Theta(\sigma)$ onto $\tensor_{i=1}^sS^{m_i}({U_{b_i}^{t_i}}^*)$ via $\varpi_1\tensor\cdots\tensor\varpi_s$ we arrive, under the natural identification of $\tensor_{i=1}^sS^{m_i}({U_{b_i}^{t_i}}^*)$ with $S^{m_1,\ldots,m_s}(W^*)$, at  the invariants  $\varphi_\sigma$ defined above. Since $S(W^*)$ is the direct sum  $\oplus_{r\in\mathbb{Z}_{\geq 0}}S^r(W^*)$, each summand of which in turn is a direct sum of $S^{m_1,\ldots, m_s}(W^*)$ as $(m_1,\ldots,m_s)$ varies over $s$-tuples of non-negative integers such that $\sum_im_i=r$, we get the required result.
\end{proof}
\renewcommand{\o}{\textup{o}}
\subsection{Restitution map and the polynomial ring on $W_\o$}\label{ss:polynomials}
For a $G$-graded vector space $V$ let $M:=V\tensor_{\mathbb C}\grass$  and $M_o$ denote the graded component in $M$ corresponding to the identity element $o\in G$. In this section we define the polynomial ring over $M_\circ$ and the `restitution map' from the space of multilinear maps on $W$ to this polynomial ring. For this, let us consider the $\grass$-module of all functions from $M\to \grass$, denoted by $\mathcal{F}(M,\grass)$. Let $F^r:T^r(M^*)\to \mathcal{F}(M,\grass)$ be given by $F^r({\boldsymbol{\alpha}})(v)=ev({\boldsymbol{\alpha}}\tensor v\tensor\cdots\tensor v)$ where $ev:T^{r}(M^*)\tensor T^r(M)\to \grass$ is as defined in \S\ref{ss:generalities}.

The symmetric group $S_r$ acts on $T^r(M)$ as described in \S\ref{action}. We then have the following analogue of  \cite[Lemma~3.11]{LZ} 
\begin{lemma}
    For $\sigma=(i\ \ i+1)\in S_r$, $F^r(\sigma.\boldsymbol{\alpha})(v)=\epsilon(|\alpha_i|,|\alpha_{i+1}|)\epsilon(|v|,|\alpha_i|-|\alpha_{i+1}|)F^r(\boldsymbol{\alpha})(v)$. In particular, $F^r(\boldsymbol{\alpha})(v)=0$ for $\boldsymbol{\alpha}\in I(M^*)$ and $v\in M_\circ$.
\end{lemma}
\begin{proof}
    The lemma follows from a simple calculations using the identity \S\ref{def:epsilon}(1)  satisfied by $\epsilon$.
\end{proof}
\noindent
Let $\mathcal{F}(M_\o,\grass)$ be the $\grass$-module of all functions from $M_\o\to \grass$. 
 Let $F^\bullet:=\oplus_{r=0}^\infty F^r:T(M^*)\to \mathcal{F}(M_\o,\grass)$. As a consequence of the above lemma,  we note that this map factors through $S(M^*)$. We continue to denote the induced map from $S(M^*)\to \mathcal{F}(M_\o,\grass)$ also by $F^\bullet$, and its restriction to $S^r(M^*)$ as $F^r$. 
 
 The next result allows us to define polynomials on $M_o$ via the restitution map.  For the purpose of the proof, we fix a basis $v_1,\ldots ,v_k$ of $V$ ordered such that $|v_i|\in G$ comes before $|v_j|\in G$ whenever $i<j$. (Here the elements of $G$ are ordered as given in the beginning of \S\ref{mixed}). Let $\phi_1,\ldots,\phi_k$ be dual to the above basis. Denote the bases of $M$ and $M^*$ corresponding to the above bases also by the same notation.
 Let $\mathcal{P}^r(M_\o)$ be the image of $F^r$ in $\mathcal{F}(M_\o,\grass)$. The space of polynomials on $M_\o$ is the image of $F^{\bullet}:=\oplus_{r=0}^\infty F^r:S(M^*)\to \mathcal{F}(M_\o,\grass)$, denoted as $\mathcal{P}(M_\o)$. Then the following proposition is a $G$-graded analogue of (\cite[Prop~3.14]{LZ} ).
 
\begin{proposition}
    The map $F^{\bullet}$ is an isomorphism from $S(M^*)\to \mathcal{P}(M_\o)$. 
\end{proposition}
\begin{proof}
     The surjectivity of the map $F^{\bullet}: S(M^*)\to \mathcal{P}(M_\o)$ is just a consequence of the definition of $ \mathcal{P}(M_\o)$. To obtain the injectivity, we show the injectivity of each $F^r$. For this we consider an $f\in \ker F^r$. The symmetric algebra $S(M^*)$ has a basis given by monomials in $\phi_i$; the monomials whose total degree is $r$ form a basis of $S^r(M^*)$. Expressing $f$ in terms of this basis, we have $$f=\sum_{r_1+\cdots+r_{k}=r}\lambda_{r_1,\ldots,r_{k}}\phi_1^{r_1}\cdots\phi_{k}^{r_{k}}.$$
     Let $v=\sum_ib_i\lambda_i\in M_\o$ for some $\lambda_i\in\grass$. Since $v\in M_\o$, $|b_i|=-|v_i|$ for all $i$. Evaluating $F^r(f)$ at $v$, we get\[F^r(f)(v)=\sum_{r_1+\cdots +r_{k}=r}\lambda_{r_1,\ldots,r_{k}}\lambda_1^{r_1}\cdots\lambda_{k}^{r_{k}}.\]Choose $N>0$ such that $\lambda_{r_1,\ldots,r_{k}}\in \grass(N)$ for all indices $r_1,\ldots,r_{k}$ such that $\sum_i r_i=r$. We now inductively choose $\lambda_j$ for $j=1,\ldots,k$ such that $\lambda_j\in (\grass(N_j)\setminus\grass (N_{j-1}))\cap X_{-|b_j|}$ for some suitable $N_j>N_{j-1}$; set $N_0=N$. For this choice of scalars $\lambda_i$, noting that the individual terms in the summation are distinct basis vectors of $\grass$, we deduce that $\lambda_{r_1,\ldots,r_{k}}=0$.    
\end{proof}
The space $\mathcal{P}^r(M_\o)$ is called the space of {\it polynomials of degree $r$ on $M_\o$}.  For $f\in S^r(M^*)$ and $w\in M_\o$, we note that $F^r(f)(w)=ev(\boldsymbol{f}\tensor w^{\otimes r})$ where $w^{\otimes r}=w \otimes w \otimes \cdots \otimes w$ ($r$-times) and $\boldsymbol{f}\in T^r(M^*)$ such that $\varpi_r(\boldsymbol{f})=f$; (recall, $\varpi_r: T^r(M^*)\to S^r(M^*)$ is the natural quotient map).

\subsection{Polynomial invariants of $W_\o$}
Let $W$ be the mixed tensor space as defined in the begining of this section. Then the graded component of $W$ corresponding to the identity element $\o\in G$ is $W_\o=\oplus_{i=1}^s({U_{b_i}^{t_i}})_\o$. As described above, we obtain the restitution map  $F^r:S^r(W^*)\to \mathcal{F}(W_\o,\grass)$.

For a tuple $(m_1,\ldots,m_s)$ such that $\sum_{i=1}^sm_i=r$, 
  let $T_\sigma\in (\tensor_{i=1}^s(U_{b_i}^{t_i})^{\tensor m_i})^*$ be the linear map corresponding to $\sigma\in S_N$ obtained in \S\ref{ss:pic}. The graded picture invariants defined in \S\ref{ss:pic} are the images in $S^{m_1,\ldots, m_s}(W^*)$ of these linear maps $T_\sigma$, $\sigma\in S_N$. Then we have,
\begin{proposition}\label{p:restitution2}
    The graded picture invariants $\varphi_\sigma\in S^{(m_1,\ldots,m_s)}(W^*)$ maps under restitution $F^r$ to the element of $\mathcal{F}(W_\o,\grass)$ given by $\boldsymbol{u}\mapsto T_\sigma(u_1^{\otimes m_1} \tensor\cdots\tensor u_s^{\otimes m_s})$ where $\boldsymbol{u}=(u_1\ldots,u_s)\in W_\o$.
\end{proposition}
\begin{proof}
    As noted earlier in \S\ref{ss:isomorphisms}, the isomorphism  $\iota:T^r(W^*)\to (T^r(W))^*$ is obtained from the non-degenerate pairing $T^r(W^*)\tensor T^r(W)\to\Lambda_{\epsilon}$ given by the evaluation map. Therefore, for any $\phi\in T^r(W^*)$ we get a linear map $\iota(\phi)$ on $T^r(W)$ and the evaluation $\iota(\phi)(w_1\tensor\ldots\tensor w_r)$ is given by $ev(\phi\tensor w_1\tensor\ldots\tensor w_r)$. In particular, $F^r(\phi)(\boldsymbol{u})=\iota(\phi)(\boldsymbol{u}\tensor\ldots\tensor \boldsymbol{u})$.
    
    The linear maps on $U_{b_i}^{t_i}$ are regarded as linear maps on $W$ via the projection $p_i: W\to U_{b_i}^{t_i}$. Denote the induced map on the dual spaces as $p_i^*:{U_{b_i}^{t_i}}^*\to W^*$.
    For a tuple $(m_1,\ldots,m_s)$ such that $\sum_{i=1}^sm_i=r$, $T^{m_1}({U_{b_1}^{t_1}}^*)\tensor\ldots\tensor T^{m_s}({U_{b_s}^{t_s}}^*)$ is a subspace of $T^r(W^*)$, via $\tensor_i {p_i^*}^{\tensor m_i}$. Similarly, $\tensor_{i=1}^s(U_{b_i}^{t_i})^{\tensor m_i}$ is a direct summand of the tensor space $T^r(W)$ so the projection   map $\textup{pr}: T^r(W)\twoheadrightarrow \left(\tensor_{i=1}^s(U_{b_i}^{t_i})^{\tensor m_i}\right) $ induces an injective map $\left(\tensor_{i=1}^s(U_{b_i}^{t_i})^{\tensor m_i}\right)^*\hookrightarrow {T^r(W)}^*$
given by  $\psi\mapsto [(w_1\tensor \ldots \tensor w_r)\mapsto\psi\circ\textup{pr}(w_1\tensor \ldots \tensor w_r)]$ for $\psi\in \left(\tensor_{i=1}^s(U_{b_i}^{t_i})^{\tensor m_i}\right)^* $ and $w_1\tensor \ldots \tensor w_r\in T^r(W)$. The non-degenerate pairing above restricted to the subspace $T^{m_1}({U_{b_1}^{t_1}}^*)\tensor\ldots\tensor T^{m_s}({U_{b_s}^{t_s}}^*)\tensor T^{m_1}(U_{b_1}^{t_1})\tensor\ldots\tensor T^{m_s}(U_{b_s}^{t_s}) $ via the above described maps, is a non-degenerate pairing and induces the isomorphism in (\ref{eq:iso}). For $r=\sum_i m_i$, $\boldsymbol{f_1}\tensor\cdots\tensor\boldsymbol{f_r}\in \tensor_{i=1}^s({U_{b_i}^{t_i}}^*)^{\tensor m_i}$ and $\boldsymbol{u}\in W_\o$, we have \begin{equation}
    \iota\circ (\tensor_{i}{p_i^*}^{\tensor m_i})(\boldsymbol{f_1}\tensor\cdots\tensor\boldsymbol{f_r})(\boldsymbol{u}\tensor\ldots\tensor \boldsymbol{u})=ev(p_1^*(\boldsymbol{f_1})\tensor \cdots\tensor p_s^*(\boldsymbol{f_r})\tensor \boldsymbol{u}\tensor\cdots\tensor\boldsymbol{u}).\end{equation} 
(Here the scaling factor involved is $1$ since $\boldsymbol{u}\in W_\o$.)  
On the other hand, the isomorphism in (\ref{eq:iso})  takes
 $\boldsymbol{f_1}\tensor\cdots\tensor\boldsymbol{f_r}$ to the linear map on $T^r(W)$ given by $w_1\tensor \cdots\tensor w_r \mapsto  ev(\boldsymbol{f_1}\tensor\cdots\tensor\boldsymbol{f_r}\tensor pr(w_1\tensor \cdots\tensor w_r))$. This linear map on  $\tensor_{i=1}^s((U_{b_i}^{t_i})^{\tensor m_i}$ also is denoted by $\iota(\boldsymbol{f_1}\tensor\cdots\tensor\boldsymbol{f_r})$.  When $w_i=\boldsymbol{u}\in W_\o$ for all $i$, the latter equals $ev(\boldsymbol{f_1}\tensor\cdots\tensor\boldsymbol{f_r}\tensor u_1^{\otimes m_1} \tensor\cdots\tensor u_s^{\otimes m_s})$. 
 This in turn evaluates to the right hand side of the above equation.
 

In the construction of $\varphi_\sigma$ in Theorem~\ref{t:pictures}, let $T_\sigma\in \left(\tensor_{i=1}^s((U_{b_i}^{t_i})^{\tensor m_i}\right)^*$ maps to $\varphi_\sigma\in S^r(W^*)$. Let $\boldsymbol{\varphi_\sigma}\in T^{m_1}({U_{b_1}^{t_1}}^*)\tensor\ldots\tensor T^{m_s}({U_{b_s}^{t_s}}^*)$ be such that $\iota(\boldsymbol{\varphi_\sigma})=T_\sigma$ and $\varpi_r(\tensor_i{p_i^*}^{\tensor m_i}(\boldsymbol{\varphi_\sigma}))=\varphi_\sigma$. By the above discussion, we have $F^r(\varphi_\sigma)(\boldsymbol{u})=ev(\tensor_i{p_i^*}^{\tensor m_i}(\boldsymbol{\varphi_\sigma})\tensor \boldsymbol{u}\tensor\ldots\tensor \boldsymbol{u})=ev (\boldsymbol{\varphi_\sigma}\tensor u_1^{\otimes m_1} \tensor\cdots\tensor u_s^{\otimes m_s})$ where $\boldsymbol{u}=(u_1\ldots,u_s)\in W_\o$. As $\iota(\boldsymbol{\varphi_\sigma})=T_\sigma$, the latter is $T_\sigma( u_1^{\otimes m_1} \tensor\cdots\tensor u_s^{\otimes m_s})$,  as required.
\end{proof}

\subsubsection{Graded picture invariants in terms of traces}
We now restrict to the case when $W=(U_1^1)^s$. Under the identification $U_1^1\cong End_{\grass} (U)$, we define a product operation on $U\tensor U^*$ as $(v\tensor \alpha). (w\tensor \beta)=v\alpha(w)\tensor\beta$ making the identification an isomorphism of $G$-graded algebras. Consider the trace function on $U_1^1$ given by $tr(v\tensor\alpha)=\epsilon(|v|,|\alpha|)\alpha(v)$. With this notation, one may define the trace monomial $tr_\sigma$ for a permutation $\sigma\in S_N$ as $tr_\sigma(v_{1}\tensor \phi_1\tensor\cdots \tensor v_{N}\tensor \phi_{N})=tr(v_{{i_1}}\tensor \phi_{{i_1}}.v_{{i_2}}\tensor \phi_{{i_2}}.\cdots .v_{{i_r}}\tensor \phi_{{i_r}})tr(v_{{j_1}}\tensor \phi_{{j_1}}.v_{{j_2}}\tensor \phi_{{j_2}}.\cdots .v_{{j_t}}\tensor \phi_{{j_t}})\cdots$
where $\sigma^{-1}=(i_1 ~i_2~ \ldots i_r)(j_1~j_2\cdots j_t)\cdots$. This definition is dependent on the permutation $\sigma$ and its cycle decomposition as well. However, if we restrict to $v_{1}\tensor \phi_1\tensor\cdots \tensor v_{N}\tensor \phi_{N}$ coming from $W_\circ$ then the definition is independent of the cycle decomposition of the permutation.  The proof of the following is analogous to that of Lemma~3.8 of \cite{PS}:
\begin{lemma}\label{l:trace} (see \cite[Lemma~4.3]{Berele2})
For a $\sigma\in S_N$ such that $\sigma^{-1}=(i_1 ~i_2~ \ldots i_r)(j_1~j_2\cdots j_t)\cdots\in S_N$, the $\gl_\epsilon(U)$-invariant map $T_\sigma$ (as in Proposition~\ref{p:restitution2}) corresponds to the trace monomials $tr_\sigma$ up to a scalar. Further, both the maps agree when restricted to the degree $0$ part, $((U_1^1)_\o)^{\tensor N}$.
\end{lemma}\hfill{$\Box$}\\
 The invariants in $\mathcal{P}(W_\o)$ for the induced action of $\gl_{\epsilon}(U)$ such that the isomorphism is $\gl_{\epsilon}(U)$-equivariant are called the invariant polynomials on $W_\o$. We recover Theorem~5.6 of \cite{Berele} as follows.
 
\newcommand{\str}{\textup{str}}
\begin{theorem}\label{t:traces}
    The invariant polynomials for the simultaneous action of $\gl_{\epsilon}(U)$ on $\oplus_{i=1}^s(U_1^1)_o$ is spanned by the trace monomials $\tr_\sigma$ given by \[\tr_\sigma(A_1,\ldots, A_s):= \tr(A_{f(i_1)}\cdots A_{f(i_r)})\tr(A_{f(i_{r+1})}\cdots A_{f(i_t}))\cdots\] where $A_1,\ldots, A_s\in\oplus_{i=1}^sU_1^1$, $\sigma=(i_1 \ldots i_r)(i_{r+1} \ldots i_t)\cdots\in S_n$ and a map $f:\{1,\ldots, n\}\to \{1,\ldots, s\}$ as $n$ varies over $\N$.
\end{theorem}
\begin{proof}
    The invariants in $\mathcal{P}(W_0)$ is the image of $S(W^*)^{\gl_{\epsilon}(U)}$ which in turn is spanned by $\varphi_\sigma$, by Theorem~\ref{t:pictures}. Proposition~\ref{p:restitution2} followed by Lemma~\ref{l:trace} then gives the required result.
\end{proof}

\vspace{1cm}

{\bf Declaration of competing interest} 

The authors declare that they have no known competing financial interests or personal relationships that could have appeared to influence the work reported in this paper.

{\bf Data availability}

No data was used for the research described in the article.



\end{document}